\documentclass[11pt]{article}

\usepackage{amssymb}
\usepackage{amsmath}
\usepackage{graphicx}
\usepackage{chemarr}

\newcommand{\be}[1]{\begin{equation}\label{#1}}
\newcommand{\ee}{\end{equation}}
\newcommand{\beqn}{\begin{eqnarray*}}
\newcommand{\eeqn}{\end{eqnarray*}}
\newcommand{\ba}{\begin{align}}
\newcommand{\ea}{\end{align}}
\newcommand{\oneh}{\frac{1}{2}}
\newcommand{\oneq}{\frac{1}{4}}
\newcommand{\ones}{\frac{1}{6}}
\newcommand{\mypmatrix}[1]{\left(\begin{array}{cccccccccccc}#1\end{array}\right)}
\newcommand{\arrowchem}[1]{\xrightarrow{#1}}
\newcommand{\arrowschem}[2]{\xrightleftharpoons[#1]{#2}}
\newcommand{\edo}{\end{document}}
\newcommand{\R}{{\mathbb R}}  
\newcommand{\Z}{{\mathbb Z}}  

\newcommand{\halmos}{\rule{1ex}{1.4ex}}
\newcommand{\qed}{\hfill \halmos} 

\newcommand{\atopnew}[2]{\genfrac{}{}{0pt}{}{#1}{#2}}
\newcommand{\Prob}{\mbox{Pr}\,}

\title{A Symbolic Computation Approach to a Problem\\
Involving Multivariate Poisson Distributions%
\footnote{%
Accompanied by Maple package {MVPoisson}
downloadable from 
\hfill\break
{\tt http://www.math.rutgers.edu/\~{}zeilberg/mamarim/mamarimhtml/mvp.html} }}

\author{Eduardo D.\ Sontag and Doron Zeilberger\\
Department of Mathematics\\
Rutgers University\\
Hill Center-Busch Campus\\
110 Frelinghuysen Rd., Piscataway, NJ 08854-8019, USA.\\
{\tt [sontag,zeilberg]  at math dot rutgers dot edu}}

\begin{document}

\maketitle

\begin{center}
{\bf Abstract}
\end{center}

Multivariate Poisson random variables subject to linear integer constraints
arise in several application areas, such as queuing and biomolecular networks.
This note shows how to compute conditional statistics in this context,
by employing WF Theory and associated algorithms.  A symbolic computation
package has been developed and is made freely available.  A discussion of
motivating biomolecular problems is also provided.

\section{Introduction}

In application areas such as queuing and biomolecular networks, one is often
interested in the study of independent Poisson random variables subject to side
information represented by linear integer constraints.  We show how to reduce
the computation of conditional statistics for this problem to the evaluation
of coefficients of generating functions. These
coefficients can, in turn, be computed using Wilf-Zeilberger (WZ) theory.  We
discuss this reduction, and make available a symbolic computation package
developed for that purpose.

We next provide a formulation of the problem, and briefly indicate its
motivations.  In Section~\ref{sec:genfun}, we explain the reduction to
exponential type generating functions, and in Section~\ref{sec:recurrences} we
discuss the fact that recurrences can be obtained for their coefficients.
Section~\ref{sec:tworow} discusses the special case of just two side
constraints, which is considerably simpler.
Section~\ref{sec:package} illustrates the use of the symbolic package
through a number of examples, all of which arise from the biomolecular
networks discussed in Section~\ref{sec:chem}.
An Appendix includes a proof of the basic representation theorem which
enables application of this techniques to certain reaction networks.

Suppose that we have $n$ independent Poisson random variables, $X_j$ 
($j=1 \ldots  n$), with parameters
$\lambda _j$ respectively. In other words
\begin{equation}
\label{eq:P1}
\Prob(X_1=k_1, X_2=k_2, \ldots  , X_n=k_n)=
e^{-(\lambda _1 + \ldots  + \lambda _n)}
\frac{\lambda _1^{k_1}}{k_1!}\frac{\lambda _2^{k_2}}{k_2!}
\ldots   \frac{\lambda _n^{k_n}}{k_n!} \quad .
\end{equation}

Suppose that we can't observe the $X_j$'s directly, 
but only a certain number, $m$, of linear combinations
of them:
\[
Y_i=\sum_{j=1}^{n} a_{ij} X_j \quad , \quad (i=1 , \ldots  , m )\quad ,
\]
where $(a_{ij})$ is a certain $m \times n$ matrix with 
non-negative coefficients.

We are interested in the following questions:

\begin{enumerate}
\item
Can one compute (fast!{}), for any given vector $(b_1, \ldots , b_m)$ 
(possibly with large coordinates), the
probability
\[
F(b_1, \ldots , b_m):= \Prob(Y_1=b_1, \ldots , Y_m=b_m) \quad .
\]
\item
Can one compute (fast!{}), for any given vector $(b_1, \ldots , b_m)$, 
(possibly with large coordinates) the
conditional expectation
\[
G_j(b_1, \ldots , b_m):= 
E[X_j \,\, \bigl \vert \, \, Y_1=b_1, \ldots  , Y_m=b_m] \quad , \quad 
(1 \leq j \leq n) .
\]
\item
More generally, can one compute (fast!{}), 
the higher moments
\[
G_j^{(r)} (b_1, \ldots , b_m):= 
E[X_j^r \,\, \bigl \vert \,\, Y_1=b_1, \ldots  , Y_m=b_m]
\quad , \quad (r \geq 2) \quad ,
\]
that would immediately allow us to compute the moments about the mean.
Can we compute (fast!{}) mixed moments, in particular the covariances?
\end{enumerate}

For example, suppose that $X_i$ is Poisson with parameter $\lambda _i$, $i=1,2$,
$X_1$ and $X_2$ are independent, and $A=(1 \; 1)$.
Thus, $Y=X_1+X_2$ is Poisson with parameter $\lambda _1+\lambda _2$.
Fix a non-negative integer $b$.
The probability that $X_1=k$ given that $Y=X_1+X_2=b$
is:
\[
e^{-(\lambda _1+\lambda _2)} \frac{\lambda _1^k}{k!}\frac{\lambda _2^{b-k}}{(b-k)!}
\;{\Bigg/} \;
e^{-(\lambda _1+\lambda _2)} \frac{(\lambda _1+\lambda _2)^b}{b!}
\]
which equals
\[
{b \choose k} \,p^k(1-p)^{b-k}
\]
with $p = \frac{\lambda _1}{\lambda _1+\lambda _2}$.
It follows that
$(X_1|Y=b)$ is a binomial random variable $B(b,p)$, and similarly $(X_2|Y=b)$
is a binomial random variable $B(b,1-p)$.
Statistics for binomial variables (means, variances, and all moments) are
well-known and easy to compute.  On the other hand, for more complicated
linear constraints, and especially if more than one such constraint is
imposed, statistics become considerably harder to obtain.

A simple example of where this type of problem might arise is as follows.
Suppose that the random variables $X_i$ count the number of calls placed,
during a typical 
time period, to an international service center and originating from a
specific country or geographical area and in a specific customer language.
For example, 
$X_1$ may represent the number of English-speaking callers from the USA, 
$X_2$ the number of Spanish-speaking callers from the USA, 
$X_3$ the number of English-speaking callers from Latin America, 
$X_4$ the number of Spanish-speaking callers from Latin America,
$X_5$ the number of English-speaking callers from the UK, 
and
$X_6$ the number of Spanish-speaking callers from the UK.
It is natural to assume that each of the random variables is
Poisson-distributed.
Now, suppose that we want to know what are the statistics of
the variable $X_1$, for example, the variance in the number of
English-speaking callers from the USA, subject to the additional information
that the total number of Spanish-language calls received was 100
and that the number of calls received from the US was 50.
This is $E[X_1|Y_1=100,Y_2=50]$ with $Y_1 = X_2+X_4+X_6$ and $Y_2=X_1+X_2$.
(More interestingly, one might have mixed information, represented by
more general linear combinations.)
We were originally motivated in this work by applications in molecular
biology; we defer to Section~\ref{sec:chem} a detailed discussion and
examples. 

\section{The generating function}
\label{sec:genfun}

Fix a matrix $A=(a_{ij})$ ($1 \leq i \leq m$,  $1 \leq j \leq n$),
once and for all. 
Let
\[
F_0(b_1, \ldots , b_m)=
\sum_{
\atopnew{k_1, \ldots , k_n \geq 0}%
{a_{11} k_1 + \ldots  + a_{1n} k_n=b_1 , \ldots  , 
a_{m1} k_1 + \ldots  + a_{mn} k_n=b_m}
}
\frac{\lambda _1^{k_1}}{k_1!}\frac{\lambda _2^{k_2}}{k_2!} \ldots   
\frac{\lambda _n^{k_n}}{k_n!}
\]
(value is zero if the sum is empty).
Thus, our focus will be on computing $F_0$, from which we can easily obtain
$F$, since
\[
F(b_1, \ldots , b_m)=
e^{-(\lambda _1 + \ldots  + \lambda _n)}F_0(b_1, \ldots , b_m) \,.
\]
Let 
$f_0$
be the (multivariable) generating function of 
$F_0$,
in other words
\[
f_0(z_1, \ldots , z_m)=\sum_{b_1 \geq 0, \ldots , b_m \geq 0} 
F_0(b_1, \ldots , b_m) z_1^{b_1} \ldots  z_m^{b_m} \quad.
\]

Our quantity of interest, $F_0(b_1, \ldots , b_m)$,
is the coefficient of $z_1^{b_1} \ldots  z_m^{b_m}$ in the multivariable
Taylor expansion about the origin of $f_0(z_1, \ldots , z_m)$.

We have:
\[
f_0(z_1, \ldots , z_m)
=\sum_{b_1 \geq 0, \ldots , b_m \geq 0} 
\left (
\sum_{
\atopnew{k_1, \ldots , k_n \geq 0}%
{a_{11} k_1 + \ldots  + a_{1n} k_n=b_1 , \ldots  , 
a_{m1} k_1 + \ldots  + a_{mn} k_n=b_m}
}
\frac{\lambda _1^{k_1}}{k_1!}\frac{\lambda _2^{k_2}}{k_2!} \ldots   
\frac{\lambda _n^{k_n}}{k_n!}
\right )
z_1^{b_1} \ldots  z_m^{b_m} \quad .
\]
By changing the order of summation, this equals
\beqn
&& \!\!\!\!\!\!\!\!\!\!\!\!\!\!\!\!\!\!\!\!
\sum_{k_1 \geq 0, \ldots , k_n \geq 0}
\frac{\lambda _1^{k_1}}{k_1!}\frac{\lambda _2^{k_2}}{k_2!} \ldots   \frac{\lambda _n^{k_n}}{k_n!}
z_1^{a_{11}k_1+ \ldots  + a_{1n} k_n} \ldots  z_m^{a_{m1} k_1 + \ldots  + a_{mn} k_n}
\\
&=&\sum_{k_1 \geq 0, \ldots , k_n \geq 0} 
\frac{(\lambda _1 z_1^{a_{11}} z_2^{a_{21}} \ldots  z_m^{a_{m1}})^{k_1}}{k_1!}
\ldots 
\frac{(\lambda _n z_1^{a_{1n}} z_2^{a_{2n}} \ldots  z_m^{a_{mn}})^{k_n}}{k_n!}
\\
&=&
\left ( \sum_{k_1 \geq 0}  \frac{(\lambda _1 z_1^{a_{11}} z_2^{a_{21}} \ldots  z_m^{a_{m1}})^{k_1}}{k_1!} \right )
\ldots 
\left ( \sum_{k_n \geq 0}  \frac{(\lambda _n z_1^{a_{1n}} z_2^{a_{2n}} \ldots  z_m^{a_{mn}})^{k_n}}{k_n!} \right )
\\
&=&\exp ( \lambda _1 z_1^{a_{11}} z_2^{a_{21}} \ldots  z_m^{a_{m1}} )
\ldots 
\exp ( \lambda _n z_1^{a_{1n}} z_2^{a_{2n}} \ldots  z_m^{a_{mn}} )
\\
&=&\exp \left ( \lambda _1 z_1^{a_{11}} z_2^{a_{21}} \ldots  z_m^{a_{m1}} + \ldots  +
\lambda _n z_1^{a_{1n}} z_2^{a_{2n}} \ldots  z_m^{a_{mn}} \right ) \quad .
\eeqn

We have just derived

{\bf Theorem 1}: 
\[
f_0(z_1, \ldots , z_m)= 
\exp \left ( 
\sum_{j=1}^{n} \lambda _j \prod_{i=1}^{m} z_i^{a_{ij}} 
\right )
\]

The conditional probability
\[
\Pr\left(X_1=k_1, X_2=k_2, \ldots  , X_n=k_n
        \,\, \bigl \vert \, \, Y_1=b_1, \ldots  , Y_m=b_m\right)
\]
is the same as the expression in (\ref{eq:P1})
divided by $F(b)$, provided that $\sum_{j=1}^{n} a_{ij} k_j=b_i$ for all $i$,
and is zero otherwise.
Recall that the $r$th factorial moment of a random variable $W$, $E[W^{(r)}]$,
is, by definition, the expectation of $W!/(W-r)!$.
We are interested in the conditional factorial moments of $X_j$ given $Y=b$,
which we will denote as
$E[X_j^{(r)}\, \bigl \vert \, Y]$.
By definition,
$E[X_j^{(r)}\, \bigl \vert \, Y]$
is the following expression divided by $F_0(b)$: 
\begin{equation}
\label{eq:moment1}
\sum_{
\atopnew{k_1, \ldots , k_n \geq 0}
{a_{11} k_1 + \ldots  + a_{1n} k_n=b_1 , \ldots  ,
a_{m1} k_1 + \ldots  + a_{mn} k_n=b_m}}
k_j(k_j-1)\ldots  (k_j-r+1)
\frac{\lambda _1^{k_1}}{k_1!}\frac{\lambda _2^{k_2}}{k_2!} \ldots 
\frac{\lambda _n^{k_n}}{k_n!} \quad .
\end{equation}
Now, expression (\ref{eq:moment1}) is the same as
the result of applying the operator
$\lambda _j^r (\frac{\partial}{\partial \lambda _j})^r$ to
$F_0(b_1, \ldots , b_m)$
\emph{when viewing the $\lambda $'s as variables and not as constants}.
On the other hand,
\[
\lambda _j^r (\frac{\partial}{\partial \lambda _j})^r f_0(z_1,\ldots ,z_m)
=
\sum_{b_1 \geq 0, \ldots , b_m \geq 0}
\lambda _j^r (\frac{\partial}{\partial \lambda _j})^r
F_0(b_1, \ldots , b_m) z_1^{b_1} \ldots  z_m^{b_m}
\]
and therefore expression (\ref{eq:moment1}) is the same as the
coefficient of $z_1^{b_1} \ldots  z_m^{b_m}$ in
$\lambda _j^r (\frac{\partial}{\partial \lambda _j})^r f_0(z_1,\ldots ,z_m)$
Since, as formal power series, we have the representation in Theorem 1.
we conclude that expression (\ref{eq:moment1}) is the same as
the coefficient of
$z_1^{b_1} \ldots  z_m^{b_m}$ in
$(\prod_{i=1}^{m} z_i^{a_{ij}})^r f_0(z)$,
which is the same as
$F(b_1-ra_{1j},b_2- ra_{2j}, \ldots , b_m-ra_{mj})$
when all $b_i-ra_{ij}\geq 0$ and zero otherwise.
In conclusion, $E[X_j^{(r)}\, \bigl \vert \, Y]$
equals 
$F_0(b_1-ra_{1j},b_2- ra_{2j}, \ldots , b_m-ra_{mj})$
divided by $F_0(b)$.
We have proved:

{\bf Theorem 2}: The conditional factorial moments
$E[X_j^{(r)}\, \bigl \vert \, Y]$ are given in terms of the
$F_0(b_1, \ldots  , b_m)$ by
\[
\lambda _j^r \, \cdot \,
\frac{F_0(b_1-ra_{1j},b_2- ra_{2j}, \ldots , b_m-ra_{mj}) }%
{F_0(b_1, \ldots , b_m)}
\]
when all $b_i-ra_{ij}\geq 0$ and zero otherwise.

So everything depends on a fast computation of the 
coefficients $F_0(b_1, \ldots  , b_m)$,  of $f_0(z_1, \ldots  , z_m)$.

By taking mixed partial derivatives, we can easily derive analogous
expressions for mixed moments, in particular, the covariances.

\section{Recurrences}
\label{sec:recurrences}

From now on, let's assume that the entries of $A$, $(a_{ij})$,
are non-negative {\bf integers}. In that case, we can write
\[
f_0(z)=\exp( Q(z_1, \ldots , z_m) ) \quad,
\]
where $Q(z_1, \ldots  , z_m)$ is the {\bf polynomial}
\[
Q(z_1, \ldots , z_m) :=
\sum_{j=1}^{n} \lambda _j \prod_{i=1}^{m} z_i^{a_{ij}} 
\quad .
\]
By Cauchy's theorem, we can express $F(b_1, \ldots , b_m)$ as
a {\bf multi-contour integral}:
\[
F(b_1, \ldots , b_m) =
\left ( \frac{1}{2 \pi i}  \right )^m
\int_{|z_1|=c} \ldots  \int_{|z_m|=c}
\frac {\exp( Q(z_1, \ldots , z_m) )} {z_1^{b_1+1} \ldots  z_m^{b_m+1} }
\, dz_1 \, \ldots  \, dz_m
\quad .
\]
By the celebrated {\bf Wilf-Zeilberger} theory (\cite{WZ}),
$F(b_1, \ldots , b_m)$ satisfies pure {\bf linear recurrences with
polynomial coefficients} in each of its arguments.
This means that for each $i$ between $1$ and $m$, there
exists a positive integer $R_i$ (the order) and polynomials
$P^{(i)}_r(b_1, \ldots , b_m)$ ($0 \leq r \leq R_i$)
such that the following
holds, for \emph{all} $(b_1, \ldots  , b_m)$:
\[
\sum_{r=0}^{R_i} P^{(i)}_r(b_1, \ldots , b_m)
F(b_1, \ldots , b_{i-1}, b_i+r, b_{i+1} , \ldots  , b_m) = 0 \quad .
\]
Once these recurrences are known, one can compute
$F(b_1, \ldots , b_m)$ in 
time linear in $b_1 + \ldots  + b_m$ and 
with constant memory allocation (one only needs to
remember, at each stage, a constant number of values).

In rare cases, the leading term of the recurrence would vanish, in which case,
we would encounter a (discrete) ``singularity'', and would not be able to go
on, since we would have to divide by $0$, but in that case one can show that
there is an alternative route, using another order of applying the
recurrences.

The {\bf Apagodu-Zeilberger}\cite{ApZ}
multi-variable extension of the Almkvist-Zeilberger\cite{AlZ}
algorithm can find such recurrences explicitly. Unfortunately,
for matrices $A$ with more than three rows, the time taken to
find such recurrences is prohibitive, but many matrices of
interest have two or three rows.  


\section{Two-Rowed matrices}
\label{sec:tworow}

If the matrix $A$ only has two rows, and its entries are only $\{0,1\}$, then
one can express $F_(b_1,b_2)$ as a \emph{single sum}. Indeed, let

%
%

\begin{itemize}
\item
$c_{01}$ be the sum of the $\lambda _j$'s for which $a_{1,j}=0$, $a_{2,j}=1$,
\item
$c_{01}$ be the sum of the $\lambda _j$'s  for which $a_{1,j}=1$, $a_{2,j}=0$,
\item
$c_{01}$ be the sum of the $\lambda _j$'s for which $a_{1,j}=1$, $a_{2,j}=1$.
\end{itemize}
Then, we have
\[
Q(z)=c_{01} z_1+c_{10}z_2+c_{11}z_1z_2 ,
\]
and so
\beqn
f_0(z_1,z_2) &=& e^{Q(z)} \;=\; \sum_{k=0}^{\infty} \frac{Q(z)^k}{k!}
\\
&=&
\sum_{\alpha \geq 0, \beta \geq 0, \gamma \geq 0} 
 \frac{ (c_{01} z_1)^{\alpha } (c_{10} z_2)^{\beta } (c_{11} z_1 z_2)^{\gamma }}%
      {\alpha ! \beta ! \gamma !}
\\
&=&\sum_{\alpha \geq 0, \beta \geq 0, \gamma \geq 0} 
 \frac{ c_{01}^{\alpha } c_{10}^{\beta } c_{11}^{\gamma } z_1^{\alpha +\gamma } z_2^{\beta +\gamma } }%
  {\alpha ! \beta ! \gamma !} \quad .
\eeqn
To get $F_0(b_1,b_2)$, we must extract the coefficient of $z_1^{b_1} z_2^{b_2}$
 which entails
$\alpha =b_1-\gamma ,\beta =b_2-\gamma $, and we have the single-sum binomial coefficient (hypergeometric) sum
(replacing $\gamma $ by $k$)
\[
F_0(b_1,b_2)=\sum_{k=0}^{\min(b_1,b_2)}  \frac{c_{11}^k c_{01}^{b_1-k} c_{10}^{b_2-k}}{k! (b_1-k)! (b_2 -k)!} \quad .
\]
Using the {\bf Zeilberger Algorithm} (\cite{Z,PWZ}) , we get the following linear recurrence:
\[
(c_{10}b_1^2+4c_{10}-2c_{10}b_2+4c_{10}b_1-c_{10}b_1b_2)F_0(b_1+2,b_2)
+
\]
\[
(-c_{11}b_1-c_{11}+b_2c_{11}+b_2c_{10}c_{01}-2b_1c_{10}c_{01}-3c_{01}c_{10})F_0(b_1+1,b_2)+
(c_{11}c_{10}+c_{01}c_{10}^2)F_0(b_1,b_2) \, = \, 0 \quad .
\]

\section{The Maple package MVPoisson}
\label{sec:package}

All this is implemented in the Maple package
{\tt MVPoisson} accompanying this article. It is
available from the webpage of this article

{\tt http://www.math.rutgers.edu/\~{}zeilberg/mamarim/mamarimhtml/mvp.html} \quad ,

where one can also find sample input and output.

We next discuss several examples of matrices $A$ and computations using
MVPoisson.  These examples, of interest in themselves, are motivated by the
biochemical networks discussed in Section~\ref{sec:chem}.

Mostly, we illustrate the use of the command ``RecsV'', which provides the
recurrences satisfied by the coefficients $F_0$, but we also show a few
examples of other commands that compute moments.

\subsection{A one-row example}

The matrix $A$ is:
\be{eqn:A1plus1}
A \;=\; (1\;\;1)\,.
\ee
As discussed in the Introduction, the conditional random variables $(X_i|Y=b)$
are binomial.
With the notations of this paper,
\[
F_0(b) \;=\; \sum_{i+j=b} \frac{\lambda _1^i}{i!} \frac{\lambda _2^j}{j!}
        \;=\; \frac{1}{b!} (\lambda _1 + \lambda _2)^b \,.
\]
This function $F_0$ clearly satisfies the following recurrence:
\[
F_0(b+1) \;=\; \frac{\lambda _1+\lambda _2}{b+1} F_0(b)
\]
with $F_0(0)=0$ and $F(1)=\lambda _1+\lambda _2$.
Indeed, for the matrix $A$ in~(\ref{eqn:A1plus1}),
the ``RecsV($A, \lambda , b$)'' command provides the following recurrence:
\[
F_0(b_1+1) \;=\; \frac{\lambda _1+\lambda _2}{1+b_1}F_0(b_1)
\]
with initial condition $F_0(1) = \lambda _1+\lambda _2$.

\subsection{A two-row example}

The matrix $A$ is:
\be{eqn:Abimolecular}
A\;=\;\mypmatrix{1 & 0 & 1\cr
0 & 1 & 1}\,.
\ee
For the matrix $A$ in~(\ref{eqn:Abimolecular}),
the ``RecsV($A, \lambda , b$)'' command provides the following two-dimensional
recurrence:
\beqn
F_0(b_1+2, b_2)
&=& 
-\;\frac{-b_2 \lambda _3 + b_1 \lambda _3 + \lambda _3 - \lambda _1 \lambda _2}{\lambda _2 (2 +b_1)}F_0(b_1+1, b_2)
\\
&+&\frac{\lambda _1 \lambda _3 }{\lambda _2 (2 + b_1)}F_0(b_1, b_2)
\eeqn
on $b_1$ and
\beqn
F_0(b_1, b_2 + 2)
&=&\frac{-\lambda _3 + b_1 \lambda _3 - b_2 \lambda _3 + \lambda _1 \lambda _2}{\lambda _1 (b_2 + 2)}F_0(b_1, b_2 + 1)
\\
&+&
\frac{\lambda _2 \lambda _3 }{\lambda _1 (b_2 + 2)}F_0(b_1, b_2)
\eeqn
on $b_2$,
with the following initial conditions:
\beqn
\mypmatrix{
F_0(1,2) & F_0(2,2)\cr
F_0(1,1) & F_0(2,1)}
&=&
\\
&&
\hspace{-100pt}
\mypmatrix{\lambda _3+\lambda _1\lambda _2 & \lambda _2\lambda _3+\oneh\lambda _2^2\lambda _1\cr
\lambda _1\lambda _3+\oneh\lambda _2\lambda _1^2 & \oneh\lambda _3^2+\lambda _2\lambda _1\lambda _3+\oneq\lambda _2^2\lambda _1^2}.
\eeqn

\subsection{Another two-row example}

The matrix $A$ is:
\be{eqn:Awater}
A = \mypmatrix{1 & 0 & 1\cr
               0 & 2 & 1}\,.
\ee

For the matrix $A$ in~(\ref{eqn:Awater}),
the ``RecsV($A, \lambda , b$)'' command provides the following two-dimensional
recurrence:
\beqn
F_0(b_1+3,b_2)
&=&
\frac{\lambda _1}{6+2b_1}F_0(b_1+2, b_2)\\
&-&
\frac{(\lambda _3^2-b_2 \lambda _3^2+2 \lambda _2 \lambda _1^2+b_1 \lambda _3^2) }{2\lambda _2(3+b_1) (2+b_1)}F_0(b_1+1,b_2)\\
&+&\frac{\lambda _1\lambda _3^2}{2\lambda _2(3+b_1)(2+b_1)}F_0(b_1, b_2)
\eeqn
on $b_1$ and
\beqn
F_0(b_1,b_2+3) &=&
\frac{\lambda _3(b_1-2-b_2)}{\lambda _1(3+b_2)}F_0(b_1, b_2+2)
\;+\;\frac{2\lambda _2}{b_2+3}F_0(b_1, b_2+1)\\
&+&\frac{2\lambda _2\lambda _3}{\lambda _1(3+b_2)}F_0(b_1, b_2)
\eeqn
on $b_2$
with the initial conditions:
\beqn
\mypmatrix{
F_0(1,3) & F_0(2,3) & F_0(3,3)\cr
F_0(1,2) & F_0(2,2) & F_0(3,2)\cr
F_0(1,1) & F_0(2,1) & F_0(3,1)}
&=&
\\
&&
\hspace{-100pt}
\mypmatrix{
\lambda _3& \lambda _1\lambda _2& \lambda _2\lambda _3\cr
\lambda _1\lambda _3& \oneh\lambda _3^2+\oneh\lambda _2\lambda _1^2&    \lambda _2\lambda _1\lambda _3\cr
\oneh\lambda _1^2\lambda _3& \oneh\lambda _1\lambda _3^2+\ones\lambda _2\lambda _1^3& \ones\lambda _3^3+\oneh\lambda _2\lambda _1^2\lambda _3
}.
\eeqn
See Figure~\ref{fig:initialconds}.

\newcommand{\dcircle}{300}
\begin{figure}[ht]
\setlength{\unitlength}{1500sp}%
\begin{center}
\begin{picture}(8127,6102)(1336,-7276)
\put(2026,-5236){\circle*{\dcircle}}
\put(3376,-5236){\circle*{\dcircle}}
\put(4726,-5236){\circle*{\dcircle}}
\put(6076,-5236){\circle{\dcircle}}
\put(7426,-5236){\circle{\dcircle}}
\put(7426,-3886){\circle{\dcircle}}
\put(6076,-3886){\circle{\dcircle}}
\put(4726,-3886){\circle*{\dcircle}}
\put(3376,-3886){\circle*{\dcircle}}
\put(2026,-3886){\circle*{\dcircle}}
\put(2026,-2536){\circle{\dcircle}}
\put(3376,-2536){\circle{\dcircle}}
\put(4726,-2536){\circle{\dcircle}}
\put(6076,-2536){\circle{\dcircle}}
\put(7426,-2536){\circle{\dcircle}}
\put(3376,-6586){\circle*{\dcircle}}
\put(4726,-6586){\circle*{\dcircle}}
\put(6076,-6586){\circle{\dcircle}}
\put(7426,-6586){\circle{\dcircle}}
\put(2026,-6586){\circle*{\dcircle}}
\put(2026,-6586){\line( 0, 1){5400}}
\multiput(8101,-5236)(117.39130,0.00000){12}{\line( 1, 0){ 58.696}}
\multiput(8101,-3886)(117.39130,0.00000){12}{\line( 1, 0){ 58.696}}
\multiput(8101,-2536)(117.39130,0.00000){12}{\line( 1, 0){ 58.696}}
\multiput(3376,-1861)(0.00000,122.72727){6}{\line( 0, 1){ 61.364}}
\multiput(4726,-1861)(0.00000,122.72727){6}{\line( 0, 1){ 61.364}}
\multiput(6076,-1861)(0.00000,122.72727){6}{\line( 0, 1){ 61.364}}
\multiput(7426,-1861)(0.00000,122.72727){6}{\line( 0, 1){ 61.364}}
\put(2026,-6586){\line( 1, 0){7425}}
\put(1936,-7261){1}%
\put(3286,-7261){2}%
\put(4636,-7216){3}%
\put(5941,-7216){4}%
\put(7291,-7261){5}%
\put(1351,-6721){1}%
\put(1351,-4021){3}%
\put(1351,-2671){4}%
\put(1351,-5371){2}%
\end{picture}%
\end{center}
\caption{Two-dimensional recursion fills-in the values of $F_0(i,j)$ at the
  locations indicated by the open circles, using the initial data given at the
  locations indicated by the filled circles.  For programming convenience,
  indices are positive integers: in the example shown, the initial
  conditions are specified for $i,j=1,2,3$.}
\label{fig:initialconds}
\end{figure}
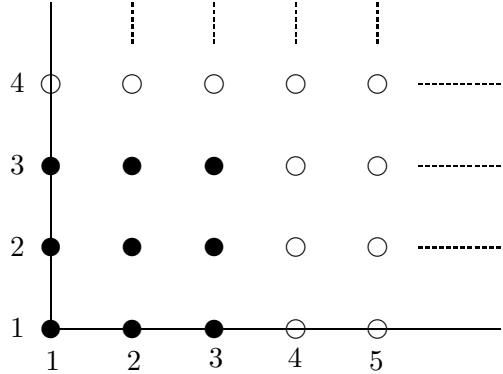

\subsection{A two-row example with five columns}

The matrix $A$ is:
\be{eqn:ARL}
A\;=\;\mypmatrix{0 & 0 & 1 & 1 & 1\cr
               1 & 1 & 0 & 1 & 1}\,.
\ee
For the matrix $A$ in~(\ref{eqn:ARL}),
the ``RecsV($A, \lambda , b$)'' command provides the following two-dimensional
recurrence:
\beqn
F_0(b_1+2, b_2)
&=&
- \;
\frac{(-b_2\lambda _5-b_2\lambda _4+\lambda _5+\lambda _4-\lambda _1\lambda _3-\lambda _2\lambda _3+b_1\lambda _5+b_1\lambda _4)}{(\lambda _1+\lambda _2)(2+b_1)}
F_0(b_1+1, b_2)\\
&+&
\frac{\lambda _3(\lambda _5+\lambda _4)}{(\lambda _1+\lambda _2)(2+b_1)}
F_0(b_1, b_2)
\eeqn
on $b_1$ and
\beqn
F_0(b_1, b_2+2) 
&=&
\frac{(-\lambda _5-\lambda _4+b_1\lambda _5+b_1\lambda _4-b_2\lambda _5-b_2\lambda _4+\lambda _2\lambda _3+\lambda _1\lambda _3)}{\lambda _3(b_2+2)}
F_0(b_1, b_2+1)
\\
&+&
\frac{(\lambda _5+\lambda _4)(\lambda _1+\lambda _2)}{\lambda _3(b_2+2)}
F_0(b_1, b_2)
\eeqn
on $b_2$, with the initial conditions:
\beqn
\mypmatrix{
F_0(1,2) & F_0(2,2)\cr
F_0(1,1) & F_0(2,1)}
&=&
\\
&&
\hspace{-120pt}
\mypmatrix{\lambda _5+\lambda _4+(\lambda _1+\lambda _2)\lambda _3 & (\lambda _5+\lambda _4)(\lambda _1+\lambda _2)+\oneh(\lambda _1+\lambda _2)^2\lambda _3\cr
\lambda _3(\lambda _5+\lambda _4)+\oneh(\lambda _1+\lambda _2)\lambda _3^2 &
\oneh(\lambda _5+\lambda _4)^2+(\lambda _1+\lambda _2)\lambda _3(\lambda _5+\lambda _4)+\oneq(\lambda _1+\lambda _2)^2\lambda _3^2}.
\eeqn

The command ``CorMf($A,\lambda ,b$)'' provides the correlation matrix for the $X_i$'s
subject to $Ax=b$ and assuming that the parameters are $\lambda $.
For the matrix $A$ considered here we obtain, for example with
$\lambda =(1,1,1,1,1)$ and $b=(5,5)$, the following result:
\[
\mypmatrix{
1.0& -.3647053019& .5636021195& -.2407443460& -.2407443460\cr
-.3647053019& 1.0& .5636021195& -.2407443460& -.2407443460\cr
.5636021195& .5636021195& 1.0& -.4271530174& -.4271530174\cr
-.2407443460& -.2407443460& -.4271530174& 1.0& -.6350805992\cr
-.2407443460& -.2407443460& -.4271530174& -.6350805992& 1.0}\,.
\]
Note the negative entry for the correlations between $X_1$ and $X_2$.
This corresponds to the fact that
$Y_2 = X_1+ X_2 + X_4 + X_5 = 5$,
so increases in $X_1$ should be expected to result in decreases in $X_2$.
Similar   interpretations apply to the other entries.

\subsection{A two-row example with six columns}

The matrix $A$ is:
\be{eqn:AGoul}
A\;=\;\mypmatrix{1 & 0 & 1 & 0& 1 & 1\cr
          0& 1 & 1 & 1 & 0 & 1}\,.
\ee
For the matrix $A$ in~(\ref{eqn:AGoul}),
the ``RecsV($A, \lambda , b$)'' command provides the following two-dimensional
recurrence:
\beqn
F_0(b_1+2, b_2)
&=&
-\;
\frac{\lambda _3+\lambda _6-\lambda _5\lambda _4-\lambda _5\lambda _2+b_1\lambda _3+b_1\lambda _6-\lambda _1\lambda _4-\lambda _1\lambda _2-b_2\lambda _3-b_2\lambda _6}{(\lambda _4+\lambda _2)(2+b_1)}
F_0(b_1+1,b_2)
\\
&+&
\frac{(\lambda _3+\lambda _6)(\lambda _5+\lambda _1)}{(\lambda _4+\lambda _2)(2+b_1)}
F_0(b_1,b_2)
\eeqn
on $b_1$, and
\beqn
F_0(b_1, b_2+2)
&=&
\frac{b_1\lambda _6+b_1\lambda _3-b_2\lambda _6-b_2\lambda _3+\lambda _1\lambda _4+\lambda _5\lambda _4+\lambda _1\lambda _2+\lambda _5\lambda _2-\lambda _6-\lambda _3}{(b_2+2)(\lambda _5+\lambda _1)}
F_0(b_1, b_2+1)
\\
&+&
(\lambda _3+\lambda _6)(\lambda _4+\lambda _2)/(b_2+2)(\lambda _5+\lambda _1)
F_0(b_1, b_2)
\eeqn
on $b_2$, 
with the initial conditions:
\beqn
F_0(1,2) &=&\lambda _3+\lambda _6+(\lambda _4+\lambda _2)(\lambda _5+\lambda _1) 
\\
F_0(2,2) &=&(\lambda _3+\lambda _6)(\lambda _4+\lambda _2)+\oneh(\lambda _4+\lambda _2)^2(\lambda _5+\lambda _1)
\\
F_0(1,1) &=&(\lambda _5+\lambda _1)(\lambda _3+\lambda _6)+\oneh(\lambda _4+\lambda _2)(\lambda _5+\lambda _1)^2 
\\
F_0(2,1) &=&\oneh(\lambda _3+\lambda _6)^2+(\lambda _4+\lambda _2)(\lambda _5+\lambda _1)(\lambda _3+\lambda _6)+\oneq(\lambda _4+\lambda _2)^2(\lambda _5+\lambda _1)^2\,.
\eeqn


\subsection{A three-row example}

The matrix $A$ is:
\be{eqn:Aenzymes}
A\;=\;
\mypmatrix{
0& 0& 1& 0& 1& 0 \cr
0& 0& 0& 1& 0& 1 \cr
1& 1& 0& 0& 1& 1 }\,.
\ee
For the matrix $A$ in~(\ref{eqn:Aenzymes}),
the ``RecsV($A, \lambda , b$)'' command provides the following two-dimensional
recurrence:
\beqn
F(b_1+2, b_2)
&=&
-\;
\frac{-b_2\lambda _6-b_2\lambda _3-\lambda _1\lambda _2-\lambda _1\lambda _4-\lambda _5\lambda _2-\lambda _5\lambda _4+b_1\lambda _6+b_1\lambda _3+\lambda _6+\lambda _3}{(\lambda _2+\lambda _4)(2+b_1)}
F(b_1+1, b_2)
\\
&+&
\frac{(\lambda _6+\lambda _3)(\lambda _1+\lambda _5)}{(\lambda _2+\lambda _4)(2+b_1)}
F(b_1, b_2)
\eeqn
on $b_1$, and
\beqn
F(b_1, b_2+2)
&=&
\frac{b_1\lambda _3+b_1\lambda _6-\lambda _3-\lambda _6+\lambda _1\lambda _2+\lambda _5\lambda _2+\lambda _1\lambda _4+\lambda _5\lambda _4-b_2\lambda _3-b_2\lambda _6}{(b_2+2)(\lambda _1+\lambda _5)}
F(b_1, b_2+1)
\\
&+&
\frac{(\lambda _6+\lambda _3)(\lambda _2+\lambda _4)}{(b_2+2)(\lambda _1+\lambda _5)}
F(b_1, b_2)
\eeqn
on $b_2$, 
with the initial conditions:
\beqn
F_0(1,2) &=&\lambda _6+\lambda _3+(\lambda _2+\lambda _4)(\lambda _1+\lambda _5)
\\
F_0(2,2) &=&(\lambda _6+\lambda _3)(\lambda _2+\lambda _4)+(1/2)(\lambda _2+\lambda _4)^2(\lambda _1+\lambda _5)
\\
F_0(1,1) &=&(\lambda _1+\lambda _5)(\lambda _6+\lambda _3)+(1/2(\lambda _2+\lambda _4))(\lambda _1+\lambda _5)^2
\\
F_0(2,1) &=&(1/2)(\lambda _6+\lambda _3)^2+(\lambda _2+\lambda _4)(\lambda _1+\lambda _5)(\lambda _6+\lambda _3)+(1/4)(\lambda _2+\lambda _4)^2(\lambda _1+\lambda _5)^2\,.
\eeqn


\subsection{A four-row example}

The $A$ matrix is:
\be{eqn:il1}
A\;=\;\mypmatrix{
1&0&0&0&1&1&0&0\cr
0&1&0&0&1&0&0&1\cr
0&0&1&0&0&1&1&0\cr
0&0&0&1&0&0&1&1}\,.
\ee
For 4-row matrices as this one, the package MVPoisson is not able to return
recurrences in a reasonable amount of time. However, one can now use the
generating functions {\it directly} to compute the relevant quantities of
interest, except that it is no longer possible to treat large inputs.

The command ``SipurD'' is used to generate averages and variances 
(``SipurD2f'' implements a more efficient algorithm specifically for
matrices with two rows).
For the matrix $A$ in~(\ref{eqn:il1}) and, for example, $\lambda =(1,1,1,1,1,1,1,1)$
we obtain that
$E[X_1|Y=b] \approx 1.897$ when $b=(10,10,10,10)$ 
and $\approx 2.813$ when $b=(20,20,20,20)$ 
(the value may be obtained to arbitrary precision),
and that the variance
of $X_1$ conditioned on this same $b$ is $\approx 1.112$ when $b=(10,10,10,10)$ 
and $\approx 1.379$ when $b=(20,20,20,20)$.
The program also guesses asymptotic formulas for these quantities
as a function of the entries of $b$, and as such is a useful tool
in research, suggesting possible general formulas that one could attempt to
prove.

\section{Biochemical applications}
\label{sec:chem}

We now explain how the problem studied here arises in the context
of systems described by chemical network theory, and in particular
chemical kinetics.
There are two fundamentally different ways to mathematically model chemical
reactions.  One of them is based on differential equations modeling, and the
other one on stochastic models.  Our problem arises from this second approach.
However, to understand its interest, it is important to first discuss the
differential equation case.
Differential equation models are useful when the number of molecules is very
large, so that a continuous approximation is appropriate.

Suppose that $n$ ``species'' interact through a network of reactions. 
The term species is used to refer to the elementary objects participating in
the interactions: in molecular biology, these are typically ions, atoms, or
molecules; in population biology and ecology, they may represent distinct
animal or plant populations, particular age groups, and so forth.
It is natural to describe such a network by a system of $n$ differential
equations which constrains the time evolution of the concentrations (or average
populations) of the various species. 
These sets of differential equations take the following general form:
\[
\dot x = \Gamma  R(x)
\]
(dot indicates time derivative)
where $x=x(t)$ is an $n$-vector of species concentrations (non-negative
real numbers) and $\Gamma $ is an $n\times m$ matrix, called the ``stoichiometry'' matrix,
whose columns describe how many units of each species are created or destroyed
by each of $m$ possible reactions.  The components of the
$m$-vector $R(x)$ quantify the reaction rates for each of the $m$ reactions, as
a function of the current concentrations as well as parameters (reaction
constants) that reflect physical and chemical information.

Chemical reactions are often described by graphs whose nodes are the
``complexes'' (the species, or combinations of species, that participate
in the reactions) and whose edges are labeled by reaction rate information.
Often, a mass-action kinetics model is used, which means that the
reaction rate is proportional to the product of the concentrations of
the reactants, and only the proportionality constant, called the kinetic
constant associated to the corresponding reaction, is displayed on an edge.
There is a systematic and simple way to map graph descriptions to
differential equations.

Some of the main results in chemical network theory were obtained by 
Horn, Jackson, and Feinberg (see ~\cite{feinberg1,feinberg2} and 
also~\cite{Tcell01} for an exposition using a somewhat different formalism). 
These results guarantee that solutions of the system of differential equations
are well-behaved (stability of equilibria, uniqueness of equilibria
modulo stoichiometric constraints), \emph{provided} that certain
structural properties are satisfied by the network.
The main such theorem is valid for what are called \emph{complex balanced}
networks.
A sufficient (though not necessary) condition for complex balancing
is that the network be ``weakly reversible'' and have ``deficiency zero''.
The deficiency is computed as $c-\ell-r$, where $c$ is the number
of complexes, $r$ is the rank of the matrix $\Gamma $, and $\ell$ is the number of
``linkage classes'' (connected components of the reaction graph).
Weak reversibility means that each connected component of the reaction graph
must be strongly connected.
We refer the reader to the citations for details on deficiency theory.
Our examples are all complex balanced.

When the numbers of molecules are very small, as is sometimes the case in
molecular biology, a discrete stochastic model may be more appropriate than a
continuous differential equation model.
Indeed, fluctuations cannot be ignored when dealing with genes (usually one
or two copies), mRNA's (in the tens), ribosomes and RNA polymerases (up to
hundreds) or certain proteins that have low concentrations.

%
%

Stochastic models fully account for the probabilistic nature of reactions.
The number of individual copies of each species at (continuous) time $t$
is viewed as a random process $X_i(t)$, $i=1,\ldots ,n$.
The Chemical Master Equation (CME), which is the differential form of the
Chapman-Kolmogorov forward equation, is a linear first-order differential
equation that describes the time evolution of the joint probability
distribution of the $X_i(t)$'s.
Often, the interest is in long-time behavior, after a transient,
that is to say in the probabilistic \emph{steady state} of the system:
the joint distribution of the random variables $X_i=X_i(\infty )$ that result in the
limit as $t\rightarrow \infty $ (provided that such a limit exists in an appropriate
technical sense).
This joint distribution is a solution of the steady state CME (ssCME), the
infinite set of linear equations obtained by setting the right-hand side of
the CME to zero. 

A very beautiful recent observation, made in~\cite{anderson08} (basically
a rewording of classical results in queuing theory in Chapter 8
of~\cite{kelly}, see also~\cite{mairesse09} for a discussion) is that
the complex balancing condition, introduced originally for deterministic
differential equation models, also guarantees that there is a solution
$\pi $ of the ssCME that is formally the joint distribution of $n$ (the
number of species) independent Poisson random variables.
More precisely, for each deterministic steady state 
$\bar x\in \R^n_{\geq 0}$
(that is, $\Gamma R({\bar x})=0$, in other words, a zero of the vector field
$\Gamma R(x)$), the vector $\pi $ defined as follows is a solution of the ssCME.
The vector $\pi $ is indexed by the $n$-dimensional lattice of non-negative
integers, $N=(N_1,\ldots ,N_n)\in \Z^n_{\geq 0}$.
We write the $N$th entry of $\pi $ as $P(N)$ (thought of as a probability, in
steady state, of the event $(X_1,X_2,\ldots ,X_n)=(N_1,\ldots ,N_n)$)
Let us write the product ${\bar x_1}^{N_1}\ldots {\bar x_n}^{N_n}$ as
``${\bar x}^N$'' and $N_1!\ldots N_n!$ as ``$N!$''.
Then, the assertion is that the vector $\pi $ whose components are
\[
P(N)\;=\; \frac{{\bar x}^N}{N!}
\]
(as well as any scalar multiple of this vector) is a solution of the ssCME.
We provide a self-contained proof of this fact in an Appendix to this paper.

However, the existence of this product form distribution does not mean that
the joint distribution of the variables $X_i$ will be independent Poisson,
because the solution of the ssCME is not, in general unique.  The lack of
uniqueness stems from conservation laws.
Because of possible conservation laws, things are a bit subtle.

As an example, suppose that two molecules of species $A$ and $B$ can
reversibly combine through a bimolecular reaction to produce a molecule of
species $C$: $A+B \leftrightarrow C$.  Let us denote the number of molecules
of species $A$, $B$, and $C$ at time $t$ by $X_i(t)$, $i=1,2,3$ respectively.
The count of $A$ molecules goes down by one every time that a reaction takes
place, at which time the count of $C$ molecules goes up by one.  Thus, the sum
of the number of $A$ molecules plus the number of $C$ molecules remains
constant: $X_1(t)+X_3(t) = b_1$.  Similarly, $X_2(t)+X_3(t)=b_2$, because the
total count of $B$ and $C$ molecules is also constant.  
This holds for all $t$, so taking limits as $t\rightarrow \infty $ (ignoring
technicalities!{}), we have that, for the steady state random variables, still
$X_1+X_3 = b_1$ and $X_2+X_3=b_2$.
Let us introduce $Y_1=X_1+X_3$ and $Y_2=X_2+X_3$.
Thus, depending on the initial conditions $b_1=X_1(0)+X_3(0)$ and
$b_2=X_2(0)+X_3(0)$, the limiting distribution will be that of $X_1$ and
$X_2$ conditioned on $Y_1=b_1$ and $Y_2=b_2$.
Once we collect this information into a matrix $A$, in this case
\[
A \;=\; \mypmatrix{1 & 0 & 1\cr
                   0 & 1 & 1}\,,
\]
we are back to the situation where we want to study the behavior of
the conditioned variables $X_i|Y_j$, where the $X_i$'s are Poisson
distributed.%
\footnote{Our discussion is incomplete from a probabilistic viewpoint, as
we have not addressed questions of uniqueness and convergence.
These questions require a careful study of irreducibility properties of the 
associated Markov chains.
We are only interested here in the computational problem of obtaining
statistics for the conditioned variables.}

The rest of this section discusses various examples.
To make the notations compatible with usage in probability theory,
we use ``$\lambda $'' for the Poisson rates (instead of ${\bar x_i}$)
and $k$ for multi-indices (instead of $N$).

\subsection{A simple reversible reaction}

Consider the following reaction:
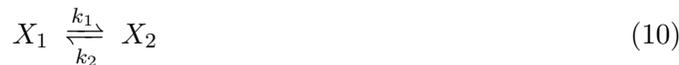
\be{eqn:1plus1}
X_1 \;\arrowschem{k_2}{k_1}\; X_2
\ee
in which one molecule of substance $X_1$ reversibly transforms to $X_2$.

This reaction system is complex-balanced, because it is weakly reversible and
it has 2 complexes, 1 strongly connected component, and rank 1, and hence
deficiency zero.

The steady states of this reaction network
are given by the solutions $\lambda =(\lambda _1,\lambda _2)$
of the equation $k_1\lambda _1=k_2\lambda _2$.  We may pick, for example,
$\lambda  = (1,k_1/k_2)$.

Every time that the forward reaction takes place, the count of molecules of
$X_1$ gets decreased by one and the count of molecules of $X_2$ gets
increased by one; the converse happens for the backward reaction.
Thus, the total number of molecules of $X_1$ and $X_2$ remains constant.
The corresponding $A$ matrix is given in~(\ref{eqn:A1plus1}).

\subsection{A bimolecular reaction}

Consider the following reaction:
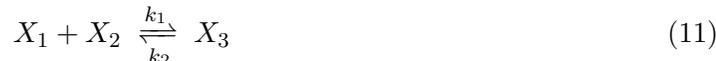
\be{eqn:bimolecular}
X_1 + X_2 \;\arrowschem{k_2}{k_1}\;  X_3
\ee
in which one molecule of $X_1$ combines reversibly with one molecule of
$X_2$ in order to produce one molecule of $X_3$.

This reaction system is complex-balanced, because it is weakly reversible and
it has 2 complexes, 1 strongly connected component, and rank 1, and hence
deficiency zero.

The steady states of this reaction network are given by the solutions
$\lambda =(\lambda _1,\lambda _2,\lambda _3)$
of the equation
\[
k_1\lambda _1\lambda _2=k_2\lambda _3 \,.
\]
We may pick, for example,
$\lambda  = (1,1,k_1/k_2)$.

Every time that the forward reaction takes place, the counts of molecules of
$X_1$ and $X_2$ each gets decreased by one and the count of molecules of $X_3$
gets increased by one; the converse happens for the backward reaction.  
Thus, the total number of molecules of $X_1$ and $X_3$ remains constant, as
does the total number of molecules of $X_2$ and $X_3$.  
The matrix $A$ is as in~(\ref{eqn:Abimolecular}).

\subsection{A more interesting bimolecular reaction}

Consider the following reaction:
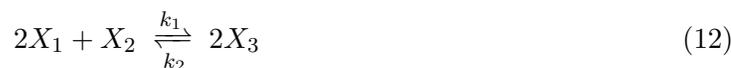
\be{eqn:water}
2X_1 + X_2 \;\arrowschem{k_2}{k_1}\; 2 X_3
\ee
which may represent, when $X_1=H_2$, $X_2 = O_2$, and $X_3 = H_2O$,
the reversible creation of a molecule of water, when two molecules of the
diatomic hydrogen gas combine with one molecule of the diatomic oxygen gas to
produce two molecules of water.  (The forward reaction produces energy, and the
reverse reaction, breaking water to form hydrogen and oxygen, requires energy,
for instance through electrolysis.  The chemical reaction formalism used here
does not account for energy production or consumption.)

This reaction system is complex-balanced, because it is weakly reversible and
it has 2 complexes, 1 strongly connected component, and rank 1, and hence
deficiency zero.

The steady states of this reaction network are given by the solutions
$\lambda =(\lambda _1,\lambda _2,\lambda _3)$
of the equation 
\[
k_1\lambda _1^2\lambda _2 = k_2\lambda _3^2 \,.
\]
We may pick, for example,
$\lambda  = (1,1,\sqrt{k_1/k_2})$.

The total sum of hydrogen and water molecules remains constant, and
for each two molecules of oxygen there is one of water produced and viceversa.
The matrix $A$ is as in~(\ref{eqn:Awater}).

\subsection{A receptor-ligand model}

\begin{figure}[ht]
\begin{center}
\setlength{\unitlength}{2600sp}%
\begin{picture}(3180,2595)(3736,-5401)
\put(4276,-3811){\vector( 0,-1){675}}
\put(4126,-4486){\vector( 0, 1){675}}
\put(6601,-4411){\vector( 0, 1){675}}
\put(6751,-3736){\vector( 0,-1){675}}
\put(5876,-3286){\vector(-1, 0){675}}
\put(5176,-3136){\vector( 1, 0){675}}
\put(5876,-5011){\vector(-1, 0){675}}
\put(5176,-4861){\vector( 1, 0){675}}
\put(6601,-5011){$C_2$}%
\put(4501,-4186){$k_{21}$}%
\put(3651,-4186){$k_{12}$}%
\put(6301,-3286){$R_2 + L$}%
\put(6926,-4186){$k_{43}$}%
\put(6126,-4186){$k_{34}$}%
\put(5401,-3586){$k_{13}$}%
\put(5401,-2986){$k_{31}$}%
\put(5401,-4711){$k_{42}$}%
\put(5401,-5311){$k_{24}$}%
\put(4126,-5086){$C_1$}%
\put(3901,-3286){$R_1 + L$}%
\end{picture}
\end{center}
\caption{A two-state receptor-ligand network}
\label{fig:RLnetwork}
\end{figure}
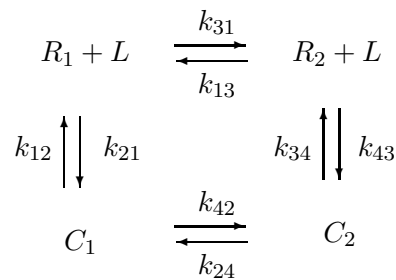

Receptor-ligand interactions play an important role in the understanding of the
biochemical mechanisms that initiate cellular signaling, and their study is
central to pharmacology.
A ``two-state'' model for such interactions studied
in~\cite{receptorligandJTB04} is shown,
pictorially, in Figure~\ref{fig:RLnetwork}.

The species participating in this reaction are:
$R_1$ and $R_2$, which represent the free receptors in an
inactive and active conformational state respectively, the free ligand $L$,
and the respective receptor-ligand complexes $C_1=R_1L$ and $C_2=R_2L$.

The steady-states $\lambda =(\lambda _1,\lambda _2,\lambda _3,\lambda _4,\lambda _5)=(R_1,R_2,L,C_1,C_2)$
of this system must satisfy the following polynomial equations:
\beqn
-(k_{21}+k_{31})R_1L +k_{12}C_1 +k_{13}R_2L&=&0
\\
-(k_{13}+k_{43})R_2L +k_{31}R_1L +k_{34}C_2 &=&0
\\
-k_{21}R_1L -k_{43}R_2L +k_{12}C_1 +k_{34}C_2&=&0
\\
-(k_{12}+k_{42})C_1 +k_{21}R_1L +k_{24}C_2 &=&0
\\
-(k_{34}+k_{24})C_2 +k_{42}C_1 +k_{43}R_2L &=&0
\eeqn
For example, when all kinetic constants are $k_i=1$ (this is not a realistic
biological choice of constants, but is picked simply for illustration), then
$\lambda =(1,1,1,1,1)$ is a steady-state.

This reaction system is complex-balanced, because it is weakly reversible and
it has 4 complexes, 1 strongly connected component, and rank 3, and hence
deficiency zero.

The conservation of $L+C_1+C_2$ (total amount of ligand) and $R_1+R_2+C_1+C_2$
(total amount of receptors) leads to the matrix in~(\ref{eqn:ARL}).

\subsection{A two-component signaling system in bacteria}

\begin{figure}[ht]
\begin{center}
\setlength{\unitlength}{1800sp}%
\begin{picture}(4530,2226)(2551,-5170)
\put(4096,-2986){\vector( 1, 0){1080}}
\put(5176,-3211){\vector(-1, 0){1125}}
\put(3601,-3661){\vector( 0,-1){900}}
\put(3376,-4561){\vector( 0, 1){900}}
\put(5176,-5011){\vector(-1, 0){1125}}
\put(7066,-4471){\vector(-3,-2){675}}
\put(6346,-4786){\vector( 3, 2){675}}
\put(6308,-3336){\vector( 3,-2){675}}
\put(2566,-5101){$R+Z$}%
\put(2566,-3230){$R+ZP$}%
\put(5416,-5101){ERP}%
\put(5416,-3230){EPR}%
\put(7066,-4066){$RP+Z$}%
\end{picture}%
\end{center}
\caption{A two-component signaling system}
\label{fig:Gnetwork}
\end{figure}
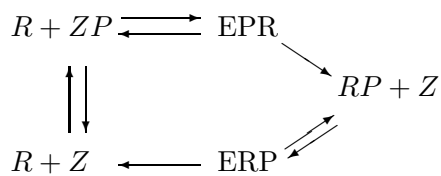

The next example is from~\cite{goulian03}.  It models the ``EnvZ/OmpR system''
in \emph{E.coli} bacteria.  This system regulates the production of certain
transport proteins (porins OmpF and OmpC) which act as pores allowing
molecules to diffuse through the cell membrane.  The system includes a kinase,
EnvZ, which phosphorylates and dephosphorylates the response regulator OmpR,
and is a particularly well-studied ``two-component signaling system'' in
bacteria. 
The model is shown, pictorially, in Figure~\ref{fig:Gnetwork}, where, for
simplicity, we omit labeling each arrow by a reaction constant.
We are using the following short-hand notations for the respective notations
in~\cite{goulian03}:
$X_1= R = \mbox{OmpR}$,
$X_2 = ZP = \mbox{EnvZ-P}$
(phosphorylated form),
$X_3 = ERP = \mbox{(EnvZ-P)OmpR}$
(complex),
$X_4 = Z = \mbox{EnvZ}$,
$X_5 = RP = \mbox{OmpR-P}$
(phosphorylated form),
and
$X_6 = EPR = \mbox{(EnvZ)OmpR-P}$
(complex).

This reaction system is complex-balanced, because it is weakly reversible and
it has 5 complexes, 1 strongly connected component, and rank 4, and hence
deficiency zero.

With all reaction constants equal to one, $\lambda =(1,1,1,1,1,1)$ is a steady state.

The system is described by six differential equations, subject to
two constraints.  These constraints reflect that the total amount of each of
OmpR and EnvZ should stay constant, respectively, and give the rows of the $A$
matrix for this example as that shown in~(\ref{eqn:AGoul}).

\subsection{A receptor antagonist model}

The paper~\cite{gnacadja07} analyzes a model involving the cytokine
Interleukin-1 (IL-1), which is produced in response to inflammatory stimuli.
The species in the model are IL-1 (denoted as L for ``ligand''), the IL-1
receptor (denoted by R), the human IL-1 receptor antagonist (denoted by A),
a decoy receptor or ``trap'' (denoted by T) which, by binding to the
ligand, helps block IL-1 signaling, and the four possible dimers RL, RA, AT,
and LT. 
The model consists of four reversible reactions:
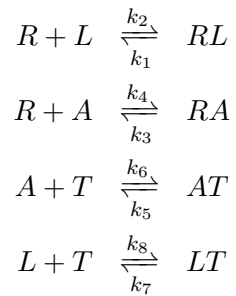
\beqn
R+L &\arrowschem{k_1}{k_2}& RL\\
R+A &\arrowschem{k_3}{k_4}& RA\\
A+T &\arrowschem{k_5}{k_6}& AT\\
L+T &\arrowschem{k_7}{k_8}& LT
\eeqn
This reaction system is complex-balanced, because it is weakly reversible and
it has 8 complexes, 4 strongly connected components, and rank 4, and hence
deficiency zero.

The total amounts of R, L, A, and T are conserved, giving rise to a matrix $A$
with 4 rows.
Ordering the states as follows:
R, L, A, T, RL, RA, AT, LT,
the resulting matrix is as in~(\ref{eqn:il1}).

\subsection{A futile-cycle example}

We now describe an example that motivates looking at a matrix
$A$ as in~(\ref{eqn:Aenzymes}).  In contrast to the previous examples,
however, this one is not complex-balanced and thus does not fit the assumptions
for the ssCME having a solution in product form.  So the interest in the
conditional statistics problem for Poisson variables is purely academic
for this particular example.
Nonetheless, it is worth seeing how such a matrix $A$ arises.

``Futile cycles'' involving phosphorylation and dephosphorylation are
ubiquitous in molecular biology (see for example~\cite{wang:JNS07} for
more discussion and references).
In such processes, an enzyme $E$ (a kinase)
catalizes the transformation of a substrate $S$ into a product $P$, passing
through one or more intermediate complexes $C$.  A different enzyme $F$ (a
phosphotase) catalizes the transformation of $P$ back into $S$, also passing
through one or mode intermediate complexes.
The simplest model (just one intermediate) for such a reaction is as follows:
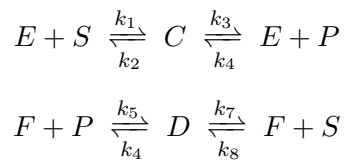
\[
E+S\;\arrowschem{k_2}{k_1}\;C\;\arrowschem{k_4}{k_3}\;E+P
\]
\[
F+P\;\arrowschem{k_4}{k_5}\;D\;\arrowschem{k_8}{k_7}\;F+S
\]
in which we used $C$ and $D$ to denote the intermediate complexes.
(Usually, the backward reactions to complex dissociation, labeled by
$k_4$ and $k_8$, are not included in the model, since they are energetically
very unfavorable.)
This system has deficiency one (6 complexes, two classes, and rank 3).  Thus,
the basic deficiency zero theory does not apply.  Interestingly, however, a
variation, ``deficiency one theory'', can be used to predict the existance of
multiple steady states for this system; see~\cite{conradi}.

There are three conservation laws, corresponding to the conservation of
total kinase, phosphotase, and substrate/product.
Ordering the variables as $S,P,E,F,C,D$, we obtain the matrix $A$ as
in~(\ref{eqn:Aenzymes}).

\appendix

\section*{Appendix}

For completeness, we show here that complex balanced reactions admit
product form equilibrium densities for their Chemical Master Equations.
The proof is basically that in~\cite{anderson08,kelly,mairesse09}.

\subsubsection*{Setup}

A \emph{chemical reaction network} is specified by:

${\cal R} = \{1,\ldots ,m\}$, the set of \emph{reactions}.

${\cal C}\subseteq \R^n_{\geq 0}$, a finite set of \emph{complexes}.

Example: if there are two reactions $1: A+B \rightarrow  C+D$ and $2: 2A+C \rightarrow  B$, then
the set ${\cal C}$ will have four elements, listing the species participating in
each: 
$(1,1,0,0)$,
$(0,0,1,1)$,
$(2,0,1,0)$,
$(0,1,0,0)$.

$S,T: \R\rightarrow {\cal C}$ are the \emph{source} and \emph{target} functions that describe
which are the reactant and product complexes, respectively.

Example: with the above reactions, 
$S(1)=(1,1,0,0)$,
$T(1)=(0,0,1,1)$,
$S(2)=(2,0,1,0)$,
$T(2)=(0,1,0,0)$.

We make the following notational convention:
for vectors $x,c\in \R^n_{\geq 0}$, $x^c:=x_1^{c_1}\ldots x_n^{c_n}$
(with $0^0=1$), and for nonnegative integer vectors
$N=(N_1,\ldots ,N_n)$, $N!:=N_1! \ldots  N_n!$.

By definition, a vector $\pi  = (P(N), N\in \Z_{\geq 0}^n)$ is a steady-state solution
of the Chemical Master Equation associated to a given reaction network
if it satisfies: 
\be{eqn:cme}
\sum_{i\in {\cal R}} P(N-T(i)+S(i)) \,A_i(N-T(i)+S(i))
\;=\;
\sum_{i\in {\cal R}} P(N) \,A_i(N)
\ee
for each $N\in \Z^n_{\geq 0}$,
where $A_i(N)$ is the $i$th ``propensity function''~\cite{gil00}: 
$A_i(N)dt$ is the probability that
reaction $i$ will occur in a small time interval $[t,t+dt]$ if the
state of the system is $N$ at time $t$.
This function is proportional to the number of ways in which the $N$ molecules
can combine to form the $i$th complex:
\[
A_i(N) = k_i \frac{N!}{(N-S(i))!} \,.
\]
The constant $k_i$ is related to the deterministic kinetic constant
of the respective reaction through division by a power of the volume in which
the reaction takes place.


A \emph{complex balanced steady state} (CBSS) with respect to the given
network and kinetic constants $k$ is an $\bar x\in \R^n_{>0}$
(which is thought of as a vector of
species concentrations) such that
the following property holds \emph{for each complex $c\in {\cal C}$:}
\be{eqn:1}
\sum_{i\in T^{-1}(c)} k_i {\bar x}^{S(i)}
\;=\;
\sum_{i\in S^{-1}(c)} k_i {\bar x}^{S(i)}
\ee
(note that one can equally well write ``${\bar x}^{c}$'' and bring this term
outside of the sum, in the right-hand side).

Complex balancing means that each ``complex'' is balanced in inflow and
outflow.  This is a Kirschoff current law (in-flux = out-flux, at each node)
when one writes a chemical network.

A counter-example to complex-balancing is this reaction network:
\[
    A \arrowchem{k_1} B , \quad\quad 2B \arrowchem{k_2} 2A
\]
(or, if one prefers reversible reactions, one may take instead
an example due to Wegsheider, $A\leftrightarrow B$ and $2A \leftrightarrow B$).
In steady state, $k_1a - 2k_2b^2 = 0$.  But complex-balancing would require that
the outflow of ``$A$'' be zero (since there are no inflows into the
``complex''  $A$), which means $k_1a = 0$, and misses the nonzero steady
states.  (One could also argue with the complex $2A$, or with
$B$, or with $2B$.) 

A \emph{complex-balanced system} is one with the property that every steady
state is complex balanced.  This concept was studied in detail by Horn and
Jackson and by Feinberg in the early 1970s.
Feinberg~\cite{feinberg1,feinberg2} showed that for a special type of system
(weakly reversible and deficiency zero), for any 
kinetic constants there is a steady state $\bar x\in \R^n_{>0}$
satisfying~(\ref{eqn:1}) (and, in fact, every other steady state will also be
complex balanced).  Let us call a system with these properties a ``Feinberg-like
system''. 

\subsubsection*{The Key Lemma}

\emph{Suppose that the reaction network is a Feinberg-like system.
Let $\bar x$ satisfy~(\ref{eqn:1}), where the $k$'s are the proportionality
factors in~(\ref{eqn:cme}).
Take any function $\alpha :{\cal C}\rightarrow \R$ on complexes.  Then:}
\be{eq:2}
\sum_{i\in {\cal R}} k_i {\bar x}^{S(i)-T(i)} \alpha (T(i)) \;=\;
\sum_{i\in {\cal R}} k_i \alpha (S(i))\,.
\ee

\noindent{\bf Proof.}
Since
\[
\sum_{i\in {\cal R}} = \sum_{c\in {\cal C}} \sum_{i\in T^{-1}(c)}
\quad\mbox{and}\quad
\sum_{i\in {\cal R}} = \sum_{c\in {\cal C}} \sum_{i\in S^{-1}(c)}
\]
it is enough to show that, \emph{for each fixed $c$}:
\[
\sum_{i\in T^{-1}(c)} k_i {\bar x}^{S(i)-c} \,\alpha (T(i))
\;=\;
\sum_{i\in S^{-1}(c)} k_i \,\alpha (S(i))\,.
\]
Since $T(i)=c$ and $S(i)=c$ in the left-hand side and right-hand side
respectively, this is the same as the CBSS condition upon multiplication by
${\bar x}^{-c}\alpha (c)$. 
\qed

\noindent{\bf Corollary.}
The vector $\Pi $ with 
\[
P(N)\;=\; \frac{{\bar x}^N}{N!}
\]
is a steady state solution of the CME.

\noindent{\bf Proof.}  Obvious using $\alpha (c) = \frac{{\bar x}^N}{(N-c)!}$.
\qed


\begin{thebibliography}{10}

\bibitem{AlZ}
G. Almkvist and D. Zeilberger,
\emph{The method of differentiating under the
integral sign}, J. Symbolic Computation {\bf 10}(1990), 571-591.
 
\bibitem{anderson08}
D.F. Anderson, G.~Craciun, and T.G. Kurtz.
\newblock Product-form stationary distributions for deficiency zero chemical
  reaction networks, 2008.
\newblock arXiv.org:0803.3042.

\bibitem{ApZ}
M. Apagodu  and D. Zeilberger,
\emph{Multi-Variable Zeilberger and Almkvist-Zeilberger Algorithms and the
Sharpening of Wilf-Zeilberger Theory },
Adv. Appl. Math. {\bf 37} (2006)(Special Regev issue), 139-152

\bibitem{goulian03}
E.~Batchelor and M.~Goulian.
\newblock {{R}obustness and the cycle of phosphorylation and dephosphorylation
  in a two-component regulatory system}.
\newblock {\em Proc. Natl. Acad. Sci. U.S.A.}, 100:691--696, 2003.

\bibitem{receptorligandJTB04}
M.~Chaves, E.D. Sontag, and R.~J. Dinerstein.
\newblock Steady-states of receptor-ligand dynamics: {A} theoretical framework.
\newblock {\em J. Theoret. Biol.}, 227(3):413--428, 2004.

\bibitem{conradi}
C.~Conradi, J.~Saez-Rodriguez, E.-D. Gilles, and J.~Raisch.
\newblock Using chemical reaction network theory to discard a kinetic mechanism
  hypothesis.
\newblock {\em IEE Proceedings Systems Biology}, 152:243 -- 248, 2005.

\bibitem{feinberg1}
M.~Feinberg.
\newblock Chemical reaction network structure and the stability of complex
  isothermal reactors - i. the deficiency zero and deficiency one theorems.
\newblock {\em Chemical Engr. Sci.}, 42:2229--2268, 1987.

\bibitem{feinberg2}
M.~Feinberg.
\newblock The existence and uniqueness of steady states for a class of chemical
  reaction networks.
\newblock {\em Archive for Rational Mechanics and Analysis}, 132:311--370,
  1995.

\bibitem{gil00}
D.~T. Gillespie.
\newblock The chemical Langevin equation.
\newblock {\em Journal of Chemical Physics}, 113(1):297--306, 2000.

\bibitem{gnacadja07}
G.~Gnacadja, A.~Shoshitaishvili, M.J. Gresser, B.~Varnum, D.~Balaban, M.~Durst,
  C.~Vezina, and Y.~Li.
\newblock Monotonicity of interleukin-1 receptor-ligand binding with respect to
  antagonist in the presence of decoy receptor.
\newblock {\em Journal of Theoretical Biology}, 244:478--488, 2007.

\bibitem{kelly}
F.~Kelly.
\newblock {\em Reversibility and Stochastic Networks}.
\newblock Wiley, New York, 1979.

\bibitem{mairesse09}
J.~Mairesse and H.-T. Nguyen.
\newblock Deficiency zero petri nets and product form, 2009.
\newblock arXiv.org:0905.3158.

\bibitem{PWZ}
M. Petkovsek, H.S. Wilf and D. Zeilberger,
\emph{A=B}, AK Peters, Wellesley, (1996).
[available on-line from the authors' websites.]

\bibitem{Tcell01}
E.D. Sontag.
\newblock Structure and stability of certain chemical networks and applications
  to the kinetic proofreading model of {T}-cell receptor signal transduction.
\newblock {\em IEEE Trans. Automat. Control}, 46(7):1028--1047, 2001.

\bibitem{wang:JNS07}
L.~Wang and E.D. Sontag.
\newblock Singularly perturbed monotone systems and an application to double
  phosphorylation cycles.
\newblock {\em J. Nonlinear Science}, 18:527--550, 2008.

\bibitem{WZ}
H.S. Wilf and D. Zeilberger,
\emph{An algorithmic proof theory for hypergeometric
(ordinary and ``q") multisum/integral identities}, Invent. Math. 
{\bf 108} (1992), 575-633.

\bibitem{Z}
D. Zeilberger,
\emph{The method of creative telescoping}, J. Symbolic Computat.
{\bf 11}, 195-204 (1991).

\end{thebibliography}

\end{document}